\documentclass[11pt]{article}
\usepackage{amsmath}
\usepackage{amsfonts}
\usepackage{amsthm}
\usepackage{verbatim}
\usepackage{hyperref}

\begin{document}

\def\N{\mathbb{N}}
\def\F{\mathbb{F}}
\def\Z{\mathbb{Z}}
\def\R{\mathbb{R}}
\def\Q{\mathbb{Q}}
\def\H{\mathcal{H}}

\parindent= 3.em \parskip=5pt

\centerline{\bf{CONTINUED FRACTIONS }} 
\centerline{\bf{ (From real numbers to formal numbers)}}
\vskip 0.5 cm
\centerline{\bf{ by A. Lasjaunias}}
\vskip 1 cm 
\noindent{\bf{1. Introduction}} 

\par The study of continued fractions is an ancient part of elementary Number Theory. It was studied by Leonhard Euler in the 18-th century. Actually, a remarkable paper from him was translated from Latin language into English and published thirty years ago \cite{Euler}.  The subject has been treated very deeply by Oskar Perron at the beginning of the 20 th century, in a famous book which has been  edited several times \cite{OP}. It can also be found in several books on Number Theory, among them a famous one  is ``An Introduction to the Theory of Numbers'' due to Hardy and Wright which was reedited many times until recently \cite{HW}.

\par In the previous cited works, continued fractions were considered for real numbers. In the 1930's, L. Carlitz initiated formal arithmetic,  replacing natural integers by polynomials over a finite field. Since then, arithmetic in function fields was developed and from the middle of last century, both number fields and function fields were studied simultaneously (see for instance Hasse's treatise on the topic \cite{H}). Our main focus will be here continued fractions in function fields, and particularly function fields over a finite field. The study of certain continued fractions in this last setting begins in the 1970's with the works of Baum and Sweet \cite{BS1,BS2}, later developed by Mills and Robbins \cite{MR}. 

\par Continued fractions are the tool to study rational approximation of numbers or functions. For general information concerning Diophantine approximation and continued fractions in power series fields, and also many references, one can see \cite{S},\cite{L} and  \cite[Chap.~9]{T}.
\vskip 0.5 cm
\noindent{\bf{2. Fields of formal power series}} 
\par Let $K$ be field and $T$ a formal indeterminate. We have the following analogy :
\begin{center}
\begin{tabular}{l c c c c c}
$\pm 1$ &$\quad$ & & $\quad$ & $K^*$ \\
$\cap $      &$\quad $ & $\quad$ & $\quad$ & $\cap$  \\
$\mathbb{Z}$ & $\quad$ & $\longleftrightarrow $ & $\quad$ & $K[T]$\\
$\cap $      &$\quad $ & $\quad$ & $\quad$ &  $\cap $\\
$\mathbb{Q}$ & $\quad$ &$\longleftrightarrow $ & $\quad$ & $K(T)$ \\
$\cap $      &$\quad $ & $\quad$ & $\quad$ & $\cap$  \\
$\mathbb{R}$ &$\quad$  &$\quad $ & $\quad$ & $K((1/T))$\\
\end{tabular}
\end{center}

We replace the number expanded in base $b$  $$x=\sum_{n\leq k
  }a_nb^n \in \mathbb{R}$$  by the formal power series  $$\alpha=\sum_{n\leq
    k} a_nT^n \in K((1/T)).$$ 
Hence a non-zero element of $K((1/T))$ can be written as
$$\alpha=\sum_{i\leq i_0}u_iT^i \quad \text{ where }\quad i\in \Z,\quad u_i\in K \quad \text{ and }\quad u_{i_0}\neq 0.$$
An ultrametric absolute value is defined over this field by $\vert \alpha \vert =\vert T \vert^{i_0}$ where $\vert T\vert$ is a fixed real number greater than 1 (not specified). Note that the field $K((1/T))$ is the completed of the field $K(T)$ for this absolute value, just as $\R$ is the completed of $\Q$ for the usual absolute value.
\newline The integer (polynomial) part of $\alpha$ is $\sum_{0\leq i\leq i_0}u_iT^i\in K[T]$ and it will be denoted $[\alpha]$ as in the classical case.
\newline We have $|x+y|\leq \max(|x|,|y|)$ for $x,y\in K((T^{-1})$ (ultrametric absolute value).
\newline If $P,Q\in K(T)$ then $|P|=|T|^{\deg(P)}$ and $|P/Q|=|T|^{\deg(P)-\deg(Q)}$.
\newline Note that, with this metric, $T^{-n}$ tends to zero when the integer $n$ tends to infinity. 
\par In the sequel, $\F_q$ is the finite field of characteristic $p$, with $q$ elements ($q=p^s$, with $s\geq 1$ integer). The field $K$ will be either $\F_q$ or $\Q$. We simply denote $\F(q)$ the field of power series in $1/T$ over $\F_q$. The importance of the fields $\F(q)$ is due to the finiteness of the base field and to the existence of the Frobenius isomorphism. Indeed, let $p$ be the characteristic of $\F_q$ and $r=p^t$ where $t\geq 0$ is an integer, then the mapping $x\rightarrow x^r$ is a field-isomorphism of $\F(q)$.

\par In the case of real numbers, a number is rational if and only if the sequence of digits is ultimately periodic. It is interesting to notice that the same is true in $\F(q)$ :
$$\alpha=\sum_{i\leq i_0}u_iT^i \in \F_q(T)\iff (u_i)_{i\leq 0} \quad \text{is ultimately periodic}.$$
A very simple rational element is obtained by the following expansion :
$$T/(T-1)=1+T^{-1}+\cdots+T^{-n}+\cdots  \quad \text{in}\quad \mathbb{F}_p(T)\quad \text{or}\quad \Q(T).$$
From this example, we can observe that the property stated above fails in $\Q((T^{-1}))$. Indeed we have by formal differentiation :
$$(T/(T-1))'=-1/(T-1)^2=-T^{-2}-2T^{-3}+\cdots-nT^{-n-1}-\cdots$$
 {\bf{Exercise 1:}} What is the power series expansion of the rational function  $(T^2/(T^2-3T+2))^5$ in $\F(5)$ ?
\newline {\bf{Exercise 2:}} Show that there is $\alpha \in \F(2)$ such that $\alpha^3=1+1/T$. Give the first digits of $\alpha$ and write an algorithm to obtain as many digits as wanted.
\vskip 0.5 cm
\noindent{\bf{3. Continuants and finite continued fractions}} 
\par Let $W=w_1,w_2,\ldots,w_n$ be a finite sequence of elements belonging to a
ring $\mathbb A$ (or simply formal variables). We set $|W|=n$ for the length of the word $W$.
We define the following operators for the word $W$.
\begin{align*}
W'&=w_2, w_3,\ldots, w_n \quad \text{ or }\quad  W'=\emptyset \quad
\text{ if } \quad \vert W\vert =1\\
W''&= w_1, w_2, \ldots, w_{n-1}\quad  \text{ or} \quad W''=\emptyset
\quad \text{if } \quad \vert W\vert =1\\
W^*&=w_n, w_{n-1},\ldots, w_1.
\end{align*}
We consider the finite continued fraction associated to $W$ to be
$$[W]=[w_1,w_2,\ldots,w_n]=w_1+\cfrac{1}{w_2+\cfrac{1}{\ddots+\cfrac{1}{w_n}}}.$$
This continued fraction is a quotient of multivariate polynomials,
usually called continuants, built from the variables $w_1,w_2,\dots,w_n$. 
$$[w_1]=w_1 \quad [w_1,w_2]=(w_1w_2+1)/w_2 $$ and $$[w_1,w_2,w_3]=(w_1w_2w_3+w_1+w_3)/(w_2w_3+1).$$
\par The continuant built on $W$ will be denoted by $\langle W\rangle$. We now recall the
definition of this sequence of multivariate polynomials. 
\newline Set $\langle \emptyset\rangle =1$. If the sequence $W$ has only one element, then we have $\langle W\rangle =W$. 
Hence, with the above notations, the continuants can be computed,
recursively on the length $\vert W\vert$, by the following
formula 
$$\langle W\rangle =w_1\langle W'\rangle +\langle (W')'\rangle  \quad \text{ for} \quad \vert W \vert \geq 2 \eqno{(1)}$$
or equivalently 
$$\langle W\rangle =\langle W''\rangle w_n+\langle (W'')''\rangle  \quad \text{ for} \quad \vert W \vert \geq 2. \eqno{(1-bis)}$$
Thus, with these notations, for any finite word $W$, the finite continued fraction $[W]$ satisfies
$$[W]=\frac{\langle W\rangle}{\langle W'\rangle}.\eqno{(2)}$$
It is easy to check that the polynomial $\langle W \rangle$ is, in a certain sense, symmetric in the variables $w_1,w_2,\ldots,w_n$ : hence we have $\langle W^*\rangle=\langle W\rangle$. This symmetry implies the classical formula 
$$[W^*]=\frac{\langle W\rangle}{\langle W''\rangle}.\eqno{(2-bis)}$$

\par In the sequel, for $n\geq 1$, we consider the finite continued fraction $z=[a_1,a_2,\cdots,a_n]$, where the $a_i$'s are positive natural integers, in the real case or polynomial in $T$ of positive degree, in the formal case. Consequently, it is clear that $z\in \Q$ in the first case or $z\in K(T)$ in the second case. For $1\leq k\leq n$, we denote $[a_k,a_{k+1},\cdots,a_n]$ by $z_k$ (hence $z_1=z$). Hence, we can write
$$z=[a_1,a_2,\cdots,a_n]=[a_1,\cdots,a_k,z_{k+1}]\quad \text{for}\quad 1\leq k\leq n-1.$$
 The $a_i$ are called the partial quotients of the expansion, while the $z_i$ are called the complete quotients or the tails of the expansion. For $1\leq k\leq n$, we also define the following elements
$$x_k=\langle a_1,\cdots,a_k\rangle \quad \text{and}\quad y_k=\langle a_2,\cdots,a_k\rangle . $$
Due to the recursive definition of the continuants, we see that both sequences $(x_k)_{1\leq k\leq n}$ and $(y_k)_{1\leq k\leq n}$ satisfy the same recurrence relation
$$K_{k+1}=a_{k+1}K_k+K_{k-1}\quad \text{for}\quad 1\leq k\leq n$$
with different initial conditions in both cases :
$$(x_0,x_1)=(1,a_1)\quad \text{or}\quad (y_0,y_1)=(0,1).$$
For $1\leq k\leq n$, according to each case, real or formal, $x_k$ and $y_k$ are positive integers or polynomial of positive degree and the rational $x_k/y_k$ is called a convergent of $z$. Moreover, we clearly have $z=x_n/y_n$.
\newline Since $z=[a_1,\cdots,a_k,z_{k+1}]$, by the property $(1-bis)$ of the continuants, we get
$$z=\frac{\langle a_1,\cdots,a_k,z_{k+1}\rangle}{\langle a_2,\cdots,a_k,z_{k+1}\rangle}=\frac{\langle a_1,\cdots,a_k\rangle z_{k+1}+\langle a_1,\cdots,a_{k-1}\rangle }{\langle a_2,\cdots,a_k\rangle z_{k+1}+\langle a_2,\cdots,a_{k-1}\rangle }.$$
This implies the following important formula :
$$\frac{x_n}{y_n}=\frac{x_kz_{k+1}+x_{k-1}}{y_kz_{k+1}+y_{k-1}}\quad \text{for}\quad 1\leq k\leq n-1.\eqno{(3)}$$
\par If $\alpha$ is an element of $\Q$ with $\alpha\geq 1$ or an element of $K(T)$ with $\vert \alpha\vert \geq \vert T\vert$ then $\alpha$ can be expanded as a finite continued fraction in a unique way. That is to say, with the above choice for the $a_i$'s in each case, we have
$$\alpha=[a_1,a_2,\cdots,a_n].$$ 
\par The only restriction to this unicity of the expansion may happen in the real case if the last partial quotient is equal to $1$. As for example :
$$[5]=[4,1]\quad \text{ or} \quad  [5,2]=[5,1,1].$$
\par In both cases, the partial quotient $a_i$ is the integral part of the complete quotients $z_i$. There are obtained by induction. So $a_1=[\alpha]$, the integer part of $\alpha$. Since $\alpha=[a_1,z_2]$, we get $z_2=1/(\alpha-a_1)$ and $a_2$ is the integer part of $z_2$, etc...This can be resumed by the formulas
 \begin{equation*}
\begin{array}{ll}
\alpha=a_1+u_1, & 0\leq u_1<1 \quad (\text{or} \quad  \vert u_1\vert <1) \\
1/u_1=z_2=a_2+u_2, & 0\leq u_2<1 \quad (\text{or} \quad  \vert u_2\vert <1) \\ 
1/u_2=z_3=a_3+u_3, & 0\leq u_3<1 \quad (\text{or} \quad  \vert u_3\vert <1) \\
\cdots \\
\end{array}
\end{equation*}

This is known as the continued fraction algorithm. This algorithm continue so long as $u_n\neq 0$. If $\alpha=u/v$ where $u,v$ are positive integers or polynomials of positive degrees,  the expansion can also be obtained by the formulas 
 \begin{equation*}
\begin{array}{ll}
u=a_1v+r_1, & 0< r_1<v \quad (\text{or} \quad 0< \vert r_1\vert <\vert v \vert) \\
v=a_2r_1+r_2, & 0< r_2<r_1 \quad(\text{or} \quad 0< \vert r_2\vert <\vert r_1 \vert) \\ 
\cdots \\
r_{n-1}=a_nr_n.
\end{array}
\end{equation*}
This is Euclid's algorithm (also use to get the greatest common divisor of the rational $u/v$) and it terminates when the last remainder is zero.
\newline Example: $\alpha=38/17$ in $\Q$.
\begin{equation*}
\begin{array}{ll}
38=2*17+4,  \\
17=4*4+1 \\ 
4=4*1.
\end{array}
\end{equation*}
and $38/17=[2,4,4]$.
\newline Example: $\alpha=T^4/(T^2+1)$ in $\F_3(T)$. 
\begin{equation*}
\begin{array}{ll}
T^4=(T^2+1)*(T^2+2)+1,  \\
T^2+2=(T^2+2)*1 \\ 
\end{array}
\end{equation*}
and $T^4/(T^2+2)=[T^2+1,T^2+2]$.
\newline (Just for curiosity : we have $T^4/(T^2+2)=[T^2-2,T^2/4+1/2]$ in $\Q(T)$.)
\par Obviously, it is important to use computer calculation to obtain these continued fraction expansions. We write below the program using Maple for two more examples first $133/37$ and then $T^{17}/(T^2-1)^6$ in $\F(17)$:   
\begin{verbatim}
a := 133; b := 37; r := 1; dfc := array(1 .. 10); 
for i to 10 while r > 0 do q := iquo(a, b); r := irem(a, b);
 dfc[i] := q; a := b; b := r end do; print(dfc);
\end{verbatim}
$$[3, 1, 1, 2, 7]$$

\begin{verbatim}
p := 17; a := T^p; b := (T^2-1)^6; r := 1; dfc := array(1 .. 15);
 for i to 15 while r <> 0 do q := `mod`(quo(a, b, T), p);
 r := `mod`(rem(a, b, T), p);
 dfc[i] := q; a := b; b := r end do; print(dfc);
\end{verbatim}
$$[T^5+6*T^3+4*T, 7*T, 11*T, 5*T, 16*T, 16*T, 14*T, 3*T, T, T, 12*T, 6*T, 10*T]$$
This last continued fraction, stated right above, is somehow surprising. The reader will understand this remark by changing $17$ into an arbitrary odd prime number $p$ and $6$ into an arbitrary integer $k$ with $1\leq k<p/2$. 

\par Concerning the restriction that we made above on the rational elements ($\alpha \geq 1$ or $\vert \alpha\vert \geq \vert T\vert$), this can easily be removed. Every rational element, in both cases, is expanded as a finite continued fraction by adding, if necessary, a first partial quotient $a_0$ (with $a_0\in \Z$ or $a_0\in K$).
\newline Examples :
$$-13/4=[-4,1,3] \quad \text{ and} \quad (T^4+T^2+1)/T^4=[1,T^2+1,T^2+2]\quad \text{ in} \quad \F(3).$$

\par To conclude this section, we state several formulas concerning the continuants. Some of them are particularly useful when dealing with continued fractions in function fields. The first one is a generalization of (1). 
For any finite sequences $A$ and $B$, of elements of $\mathbb A$, defining $A,B$ as the concatenation of sequences $A$ and $B$, we have
$$\langle A,B \rangle =\langle A\rangle \langle B \rangle +\langle A''\rangle \langle B'\rangle.\eqno{(4)}$$
If $\vert A \vert =1$ then $(4)$ simply reduces to $(1)$. The general case is obtained by induction on $\vert A \vert$.
\newline Secondly, using induction on $\vert W \vert$, we have the
following classical identity 
$$\langle W \rangle \langle (W')''\rangle-\langle W'\rangle \langle W''\rangle=(-1)^{ \vert W \vert}\quad \text{ for} \quad \vert W \vert \geq 2.\eqno{(5)}$$
Note that this implies, with the above notation for rational continued fractions in $\Q$ or $K(T)$,  
$$x_ny_{n-1}-y_nx_{n-1}=(-1)^n\quad \text{ for} \quad n\geq 1. \eqno{(6)}$$
With the above convention $x_0=1$ and $y_0=0$. We observe tat $(6)$ implies $GCD(x_n,y_n)=1$. Combining the above formulas $(4)$ and $(5)$, we also have 
$$\langle A,B \rangle \langle A'\rangle-\langle A\rangle \langle A',B\rangle=(-1)^{ \vert A \vert-1}\langle B'\rangle .\eqno{(7)}$$
Coming back to rational continued fraction, for $0\leq m <n$, we take $A=a_1,\cdots,a_m$ and $B=a_{m+1},\cdots,a_n$, then $(7)$ implies the following generalization of $(6)$
$$x_ny_m-y_nx_m=(-1)^{m-1}\langle a_{m+2},\cdots,a_n \rangle \quad \text{ for} \quad 0\leq m<n. \eqno{(8)}$$
Let $y$ be an invertible element of $\mathbb A$, then we define $y
\cdot W$ as the following sequence
$$y\cdot W = y w_1, y^{-1}w_2,\ldots, y^{(-1)^{n-1}}w_n.$$
With this notation, it is easy to check that we have $y[W]=[y \cdot
W]$ and more precisely (using induction)
$$\langle y \cdot W\rangle =\langle W\rangle  \quad \text{ if} \quad \vert W \vert \text{ is
  even}$$
and $$\langle y \cdot W\rangle =y\langle W\rangle \text{ if} \quad \vert
W \vert \text{ is odd} . $$
\par Finally, as an exercise, we give the following result (derived from (3) and (6)):
\newline {\bf{Exercise 3:}}For a finite word $W$ and two formal variables $x$ and $y$, we have 
$$[W,x]=[W]+y \quad \iff \quad y \langle W'\rangle (x\langle W'\rangle +\langle (W')''\rangle)=(-1)^{\vert W \vert-1}.$$ 
\vskip 0.5 cm
\noindent{\bf{4. The golden ratio and other quadratic elements}} 
\par Every irrational number or irrational formal power series $\alpha$ can be expanded as an infinite continued fraction: $\alpha=[a_0,a_1,a_2,\cdots,a_n,\cdots]$.
\newline As above, we assume that $\alpha \geq 1$ or $\vert \alpha\vert \geq \vert T\vert$, thus we have 
$$\alpha=[a_1,\cdots,a_n,\cdots]\quad \text{and}\quad \frac{x_n}{y_n}=[a_1,\cdots,a_n]$$ where the partial quotients are positive integers or polynomials in $T$ of positive degree and the integers (or polynomials) $x_n$ and $y_n$ are the continuants defined in the previous section. These rationals $x_n/y_n$ are called the convergents of $\alpha$. The partial quotients $a_n$ are obtained by the continued fraction algorithm (extracting the integer part of the complete quotient), which never terminates due to the irrationality of the element.
\newline A famous example: $\pi \approx 3.14159265\cdots \implies \pi=[3,7,15,1,\cdots].$
\par We will denote the $n$-th complete quotient by $\alpha_n$ ($\alpha_1=\alpha$). $\alpha$ being irrational, each complete quotient is also an infinite continued fraction, so we have 
$\alpha_{n+1}=1/(\alpha_n-a_n)$ and $1/\alpha_{n+1}$ is the fractional part of $\alpha_n$. Hence we have $\alpha=[a_1,a_2,\cdots,a_n,\alpha_{n+1}]$ and this, like formula (3) given above, implies the following formula 
$$\alpha=\frac{x_n\alpha_{n+1}+x_{n-1}}{y_n\alpha_{n+1}+y_{n-1}}\quad \text{for}\quad n\geq 1.\eqno{(9)}$$
As stated above $x_n$ and $y_n$ are the continuants $\langle a_1,\cdots,a_n\rangle$ and $\langle a_2,\cdots,a_n\rangle$ respectively, and $x_0=1$, $y_0=0$.
According to formula (9) and (6), in both cases, we can write
$$\alpha-\frac{x_n}{y_n}=\frac{x_n\alpha_{n+1}+x_{n-1}}{y_n\alpha_{n+1}+y_{n-1}}-\frac{x_n}{y_n}=\frac{(-1)^{n-1}}{y_n(y_n\alpha_{n+1}+y_{n-1})} $$
and consequently we obtain
$$\vert \alpha-x_n/y_n\vert =1/\vert y_n(y_n\alpha_{n+1}+y_{n-1})\vert.\eqno{(10)} $$
In both cases $X_n=\vert y_n(y_n\alpha_{n+1}+y_{n-1})\vert$ tends to infinity with $n$ (we always have $X_n>\vert y_n\vert ^2$ and $y_n$ or $\vert y_n\vert$ tends to infinity). This shows that the sequence of the convergents tends to $\alpha$, as it is expected. 
\par Let us consider the most simple and most famous example in $\R$:
$$\phi=[1,1,\dots,1,\dots]\quad \text{ satisfying } \quad \phi=1+1/\phi.$$
Thus we see that that $\phi$ is a solution of the following quadratic equation : $$x^2-x-1=0.$$
This equation has two real roots $(1+\sqrt{5})/2$ and $(1-\sqrt{5})/2$, the first one has integer part $1$ and it is $\phi$, the second one is negative it is $-1/\phi$. Hence we have 
$$(1-\sqrt{5})/2=-[0,1,1,\cdots,1,\cdots].$$
Note that the integer part of $(1-\sqrt{5})/2$ is $-1$. To transform this formula into a proper continued fraction, we observe that $-1/\phi=-1+1/(1+\phi)$. Hence we get
$$-1/\phi=-1+1/(2+1/\phi)=[-1,2,1,1,\dots,1,\dots].$$
Note that $\phi$ is the limit of the finite continued fraction $[1 ( n\quad \text{times})]$. If we denote $f_n=\langle 1 (n \quad \text{times}) \rangle$, we see that this continuant is the $n$-th term of the famous Fibonacci sequence $(f_n)_{n\geq 0}$ defined recursively by
$$f_0=1, \quad f_1=1 \quad \text{ and } \quad f_{n+1}=f_n+f_{n-1}.$$
 So we recover the classical formula $\phi= \lim_{n \rightarrow \infty}f_n/f_{n-1}\approx 1.61803\cdots$.
Another simple example is obtained as follows. We set $u=\sqrt{2}+1$. We have $u=2+1/u$ and consequently $u=[2,2,\dots,2,\dots]$. Since $\sqrt{2}-1=1/u$, we get $\sqrt{2}=1+1/u=[1,2,2,\dots,2,\dots].$
\par Let us now consider the analogue of $\phi$ in the formal case. The following infinite continued fraction
$$\Phi=[T,T,\cdots,T,\cdots]$$ 
can be regarded as an element of $\Q((T^{-1}))$ or of $\F(p)$ (for all primes $p$). This element satisfies $\Phi=T+1/\Phi$ and therefore is solution of the quadratic equation $x^2-Tx-1=0$. The second solution of this equation is $-1/\Phi$ and here we have
$$-1/\Phi=[0,-T,-T,\cdots,-T,\cdots].$$
(Remember: in $\F(2)$, we have $+1=-1$ and the two solutions are  $\Phi$ and $1/\Phi$.)
\newline We denote $F_n=\langle T (n \quad \text{times}) \rangle$, we see that this continuant is the $n$-th term of a sequence $(F_n)_{n\geq 0}$, generalizing the Fibonacci sequence and defined recursively by
$$F_0=1,\quad F_1=T \quad \text{ and } \quad F_{n+1}=TF_n+F_{n-1}.$$
 We can prove that $F_n=\sum_{0\leq k \leq n/2}\binom{n-k}{k}T^{n-2k}$. Since we have  $[T ( n \quad \text{times})]=F_n/F_{n-1}$, we get $\Phi= \lim_{n \rightarrow \infty}F_n/F_{n-1}$. Several other formulas, similar to classical ones for Fibonacci numbers, are known. For instance, for $k\geq 1$, we have
$$\Phi^k= \lim_{n \rightarrow \infty}F_{n+k}/F_{n} \quad \text{ or } \quad \Phi^{k+1}=\Phi F_k+F_{k-1}.$$
The importance of this element $\Phi$ is strengthened by its relation to other algebraic (non-quadratic) continued fractions in $\F(q)$, where the partial quotients are polynomial in $T$ of degree one (see   \cite[p.~400]{MR} and \cite{LY,L1}).
\newline {\bf{Exercise 4:}} What are the coefficients of $\Phi$ in the power series expansion in $\F(p)$, if $p=2$ (and if $p>2$) ?
\newline {\bf{Exercise 5:}} Let $p$ be an odd prime and $r=p^t$ with $t\geq 1$ an integer. Show that we have $F_{r-1}=(T^2+4)^{(r-1)/2}$.
\par These particular quadratic continued fractions have a special form. Let us generalize, considering:
$$\alpha=[a_1,a_2,\cdots,a_n,b_1,b_2,\cdots,b_m,b_1,b_2,\cdots,b_m,b_1,b_2,\cdots,b_m,\cdots].\eqno{(UP)}$$
The sequence of partial quotients for $\alpha$ is ultimately periodic: after a finite number of terms, we get a sequence of length $m$ repeated infinitely. Then considering
$$\alpha_{n+1}=\beta=[b_1,b_2,\cdots,b_m,b_1,b_2,\cdots,b_m,b_1,b_2,\cdots,b_m,\cdots]$$
we have $\alpha=[a_1,a_2,\cdots,a_n,\beta]$ and $\beta=[b_1,b_2,\cdots,b_m,\beta]$. applying formula (9), we get $\alpha=f(\beta)$ and $\beta=g(\beta)$, where $f$ and $g$ are linear fractional transformation with integer coefficients (rational integers or polynomial in $T$). This implies that $\beta$ is quadratic and finally that $\alpha$ is also quadratic. Hence all the expansion of the form $(UP)$ are quadratic elements.
\newline {\bf{Exercise 6:}} Let $\alpha=[5,1,2,3,1,2,3,\dots,1,2,3,\dots]\in \R$, give the quadratic equation satisfied by $\alpha$.

\par In the real case, it was established by Lagrange in the 18-th century that the converse is also true: every irrational quadratic real number has a continued fraction expansion where the partial quotients form an ultimately periodic sequence. This is also true in the formal case if $K=\F_q$. This subject is discussed deeply in general function fields by Schmidt in \cite{S}.  
\par Given a particular irrational element, the continued fraction expansion is seldom known completely. This is the case for quadratic elements. In the real case, the famous transcendental number $e$, base of the neperian logarithm, is exceptional and it has a very structured sequence of partial quotients, established by Euler (\cite[p.~311]{Euler}). Concerning other algebraic real numbers (not quadratic), not a single one is known. There is a famous question concerning these algebraic real numbers of degree $\geq 3$: is the sequence of partial quotients unbounded ? This question was put forward by Kintchine in the 1930's. Computer calculations incline to believe that the answer to this question is positive. However, this is probably one of the most important unsolved questions in Number Theory. For formal power series the situation is different. In \cite{BS1}, Baum and Sweet described a famous cubic element of $\F(2)$ having partial quotients of degrees only one or two. Later, Mills and Robbins \cite{MR} could fully describe the expansion for this element and they also gave other algebraic non-quadratic elements in different characteristic, having an explicit expansion with bounded partial quotients.  
\vskip 0.5 cm
\noindent{\bf{5. The irrationality measure}} 
\par Around 1850, Liouville studied rational approximation to algebraic real numbers. He could prove the following theorem :
\newline If $\alpha$ is a real number, algebraic of degree $d>1$ over $\Q$, then there is a positive real constant $C$ depending on $\alpha$ such that $\vert \alpha-p/q\vert >Cq^{-d}$, for all pairs $(p,q)\in \Z^2$ with $q>0$.
\par About hundred years later, Mahler \cite{M} adapted Liouville's theorem in the frame of function fields. He obtained the following :
\newline If $\alpha \in K((T^{-1}))$ is algebraic of degree $d>1$ over $K(T)$, then there is a positive real constant $C$ depending on $\alpha$ such that $\vert \alpha-P/Q\vert >C\vert Q\vert ^{-d}$, for all pairs $(P,Q)\in K[T]$ with $Q\neq 0$.
\par After a long chain of works due to Thue, Siegel, Dyson and Schneider, concerning rational approximation to algebraic real numbers, the following was proved in 1955 by Roth \cite{R}:  
\newline If $\alpha$ is an irrational algebraic real number then, for all $\epsilon >0$, there is a positive real constant $C$ depending on $\alpha$ and $\epsilon$ such that $\vert \alpha-p/q\vert >Cq^{-2-\epsilon}$, for all pairs $(p,q)\in \Z^2$ with $q>0$.
\newline This means that the algebraic real numbers are poorly approximable by rational numbers.
\par Shortly after, Roth's theorem was adapted by Uchiyama \cite{U} in certain function fields :
\newline Let $K$ be a field of characteristic zero. If $\alpha \in K((T^{-1}))$ is an irrational  algebraic then, for all $\epsilon >0$, there is a positive real constant $C$ depending on $\alpha$ and $\epsilon$ such that $\vert \alpha-P/Q\vert >C\vert Q\vert ^{-2-\epsilon}$, for all pairs $(P,Q)\in K[T]$ with $Q\neq 0$.

\par When Mahler wrote his paper on the adaptation of Liouville's theorem \cite{M}, he observed that the first improvement obtained by Thue was not valid in function fields of positive characteristic. Indeed in these fields $\F(q)$ there are algebraic elements that are very well approximable, contrarily to what happens for algebraic real numbers. Mahler gives the following example:
$$\alpha=T^{-1}+T^{-r}+T^{-r^2}+\cdots+T^{-r^n}+\cdots$$
For $n\geq 1$, we set $$U_n=T^{r^{n-1}}(T^{-1}+T^{-r}+T^{-r^2}+\cdots+T^{-r^{n-1}})\quad \text{ and }\quad V_n=T^{r^{n-1}}.$$
Then we have $$\vert \alpha-U_n/V_n\vert =\vert T^{-r^n}+T^{-r^{n+1}}+\cdots\vert =\vert V_n\vert ^{-r}.$$
But $\alpha$ is irrational and it satisfies the following algebraic equation of degree $r$: $\alpha=T^{-1}+\alpha^r$. If the exact algebraic degree of $\alpha$ were $d<r$, we would have a contradiction with Mahler's theorem, hence we have $d=r$. This shows that there cannot be an improvement (and particularly no analogue of Roth's theorem) in function fields of positive characteristic. 
\par It appears necessary to define a quantity which will measure the quality of the rational approximation of a real number or a formal power series. Let $\alpha \in \R$ or $\alpha \in K((T^{-1}))$. We define the irrationality measure (or the approximation exponent) of $\alpha$ as 
$$\nu(\alpha)=-\limsup_{\vert Q \vert \to \infty}(\log \vert \alpha -P/Q \vert/\log \vert Q \vert),$$
where $P$ and $Q$ run over polynomials in $K[T]$ with $Q\neq 0$ in the formal case and $P$ and $Q$ run over rational integers with $Q>0$ in the real case.
\newline If $\beta$ is a linear fractional transformation with integer coefficients of $\alpha$, then $\alpha$ and $\beta$ have the same irrationality measure.
\par With this notation, Liouville's theorem and Malher's adaptation say, in both cases, that we have $\nu(\alpha)\leq d$ if $\alpha$ is an algebraic element of degree $d>1$. This way, by considering elements of infinite irrationality measure, Liouville could prove the existence of transcendental real numbers. A famous example : $x=\sum_{n\geq 1}10^{-n!}$ (see \cite{HW}).
\newline We will see that for all elements in $\R$ or $K((T^{-1}))$, we have $\nu(\alpha)\geq 2$ and for all $\mu\in [2;\infty]$ there is $\alpha$ such that $\mu=\nu(\alpha)$ in both cases.
\par We shall now see that this quantity $\nu(\alpha)$ can be expressed in a simpler way using the continued fraction expansion. So in the following $\alpha$ is an irrational element defined by the infinite continued fraction $[a_1,a_2,\cdots,a_n,\cdots]$. We recall that the convergents are the best rational approximation to $\alpha$. This means that, for $n\geq 1$, we have 
$$\vert \alpha-x_n/y_n \vert \leq  \vert \alpha-P/Q \vert \quad \text{ if }\quad \vert Q\vert \leq \vert y_n\vert \quad (\text{or} \quad Q \leq y_n).$$
Consequently we obtain
$$\nu(\alpha)=-\limsup_{n}(\log \vert \alpha -x_n/y_n \vert/\log \vert y_n \vert).$$
 We shall now use formula $(10)$ of Section 4.
$$\vert \alpha-x_n/y_n\vert =1/\vert y_n(y_n\alpha_{n+1}+y_{n-1})\vert.\eqno{(10)} $$

\par In the formal case, we have $\vert y_n\alpha_{n+1}+y_{n-1}\vert =\vert y_n a_{n+1}\vert$, hence (10) implies
$$\vert \alpha-x_n/y_n \vert=\vert y_n\vert^{-2}\vert a_{n+1}\vert^{-1}=\vert y_n\vert^{-2-\deg(a_{n+1})/\deg(y_n)}.$$
Note that $\deg(y_n)=\deg(a_n)+\deg(y_{n-1})$. Hence we have $\deg(y_n)=\sum_{2\leq k\leq n}\deg(a_i)$, and we get
$$\nu(\alpha)=2+\limsup_{n\geq 1}\big ( \deg(a_{n+1})/\sum_{1\leq i\leq n}\deg(a_i)\big ).\eqno{(11)}$$
\newline {\bf{Exercise 7:}} Let $u=(u_n)_{n\geq 1}$ be a sequence of positive integers and $\alpha(u)$ be the continued fraction in $\F(p)$ such that $a_n=T^{u_n}$ for $n\geq 1$. Given a real number $x>2$ find a sequence $u$ such that $\nu(\alpha(u))=x$ (Hint: start with $x\in \N$).

\par In the real case, things are slightly different, due to the absolute value. Since $a_{n+1}<\alpha_{n+1}<a_{n+1}+1$, we can write 
$$y_n\alpha_{n+1}+y_{n-1}>y_na_{n+1}+y_{n-1}=y_{n+1}$$
and
$$y_n\alpha_{n+1}+y_{n-1}<y_n(a_{n+1}+1)+y_{n-1}=y_{n+1}+y_n<2y_{n+1}.$$
Hence (10) implies 
$$1/(2y_ny_{n+1})<\vert \alpha-x_n/y_n\vert <1/(y_ny_{n+1}).$$
This leads to 
$$\nu(\alpha)=1+\limsup_{n\geq 1}\big (\log(y_{n+1})/\log(y_n)\big ).$$
\par One can prove the following statement : If $\alpha \in \R -\Q$ and we have
$\quad \lim_n(\sum_{1\leq k\leq n} \ln(1+1/a_k))/(\sum_{1\leq k\leq n}
\ln a_k)=0$ 
then $$\nu(\alpha)=2+\limsup_n\frac{\ln a_{n+1}}{\sum_{1\leq k\leq n}
\ln a_k}.$$
From this we can get the following :
\par Let $a$ and $b$ be two real numbers greater than 1, and $\alpha\in \R$ defined by $\alpha=[a_1,a_2,\cdots,a_n,\cdots]$.
\newline $\bullet \quad a_n=[a^n]$ for $n\geq 1$ $\implies \nu(\alpha)=2$.
\newline $\bullet \quad a_n=[a^{b^n}]$ for $n\geq 1$ $\implies \nu(\alpha)=b+1$.
\newline $\bullet \quad a_n=[a^{n^n}]$ for $n\geq 1$ $\implies \nu(\alpha)=+\infty$.
\par In the sequel we will only be concerned with formal power series over a finite field (or simply called formal numbers).

\vskip 0.5 cm
\noindent{\bf{6. Hyperquadratic formal power series}} 

\par As we observed in the second section, the particularity of the fields $\F(q)$, due to the positive characteristic, is the existence of the Frobenius isomorphism. This was put forward first by Mahler in his fundamental paper on Diophantine approximation \cite{M}. In the 1970's and 1980's, Osgood, Voloch \cite{V} and de Mathan \cite{dM} have considered, on this subject, a particular class of algebraic power series.  Simultaneously, Baum and Sweet and later Mills and Robbins \cite{BS1,BS2,MR} also considered the same class of algebraic elements from the point of view of continued fraction expansion. In a later work, we have called these elements hyperquadratic power series \cite{BL}. Let us give here the following definition.
\newline Let $\F(q)$ be the field of power series over the finite field $\F_q$ of characteristic $p$. An irrational element $\alpha \in \F(q)$ is called hyperquadratic if it satisfies an algebraic equation of the following form :
$$A\alpha^{r+1}+B\alpha^r+C\alpha+D=0$$
where $A,B,C$ and $D$ are in $\F_q[T]$ and $r=p^t$ for an integer $t\geq 0$.
\newline To be more precise, we say $\alpha$ is hyperquadratic of order $t$. Also the subset in $\F(q)$ of such algebraic elements can be denoted by $\H(q)$.
\par Note that if $t=0$ then $r=1$ and $\alpha$ is simply irrational quadratic. Moreover if $\alpha$ is irrational algebraic of degree 2 or 3, then $\alpha^{p+1}$, $\alpha^p$, $\alpha$ and $1$ are linked over $\F_q(T)$ and therefore $\alpha$ is hyperquadratic of order 1. Note that the above equation can also be written under the form $\alpha=f(\alpha^r)$ where $f$ is a linear fractional transformation with integer (polynomial) coefficients. 
\newline (By iteration $\alpha=f(\alpha^r)$ implies $\alpha=g(\alpha^{r^2})$ where $g$ is another linear fractional transformation with integers coefficients. So we see that if $\alpha$ is hyperquadratic of order $t$ then it is also of order $nt$ for $n\geq 1$.)
\par First example (Mahler's example) : Here $p$ is a prime and $r=p^t$ with $t\geq 1$ integer. Let $\alpha\in \F(p)$ be defined by $$\alpha=T^{-1}+T^{-r}+T^{-r^2}+\cdots+T^{-r^n}+\cdots $$ We have $\alpha=T^{-1}+\alpha^r$. Note that if $u\in \F_p^*$ then $\alpha+u$ is also solution of the same equation. Remember, as we have seen above, that this element has algebraic degree $r$ and also irrationality measure $r$.  
\par Second example (dual of Mahler's example): Again  $p$ is a prime and $r=p^t$ with $t\geq 1$ integer. Let $\alpha\in \F(p)$ be defined by $$\alpha=[T,T^{r},T^{r^2},\cdots,T^{r^n},\cdots].$$ We have $\alpha=T+1/\alpha^r$. Note that we have
$$\deg(a_n)/\sum_{1\leq i\leq n-1}\deg(a_i)=r^{n-1}(r-1)/(r^{n-1}-1)\quad \text{ for }\quad n\geq 2$$
and consequently formula $(11)$ implies $\nu(\alpha)=2+r-1=r+1$. Hence, again the irrationality measure is maximal and equal to the algebraic degree.
\par Schmidt \cite{S}, observed that this last example can be generalized as follows. Let $l\geq1$ be an integer and $(a_1,a_2,\cdots,a_l)\in (\F_q[T])^l$. Then consider the infinite continued fraction
$$\alpha=[a_1,a_2,\cdots,a_l,a_1^r,a_2^r,\cdots,a_l^r,\cdots,a_1^{r^n},a_2^{r^n},\cdots,a_l^{r^n},\cdots]. \eqno{(S)}$$
We have $\alpha_{l+1}=\alpha^r$ and therefore, by $(9)$,  $\alpha=(x_{l}\alpha^r+x_{l-1})/(y_{l}\alpha^r+y_{l-1}).$ 
\newline Schmidt, and also Thakur independently, could prove that for any rational number $\mu >2$, there is an $\alpha$ defined by $(S)$, if $l$ and the degrees of the $a_i's$ for $1\leq i\leq l$ are well chosen, such that we have $\nu(\alpha)=\mu$.
\newline {\bf{Exercise 8:}} Let $\alpha \in \F(p)$ be such that $\alpha=[a_1,a_2,\alpha^r]$ where $a_1,a_2 \in \F_p[T]$ and $d_1=\deg(a_1)$, $d_2=\deg(a_2)$. Assuming that we have $d_1>d_2\geq 1$, compute $\nu(\alpha)$. Show that $(r-1)d_1>(r-2)(d_1+d_2)$ implies that $\alpha$ is algebraic of degree $r+1$. Same questions if $d_1<d_2<rd_1$.

\par Mills and Robbins \cite[p.~402-403]{MR} proved the following result :
\newline If $\alpha=(U\alpha^r+V)/(W\alpha^r+Z)$ and $r>1+\deg(UZ-VW)$ then $\alpha$ has a sequence of partial quotients of unbounded degrees. Particularly if  $UZ-VW\in \F_q^*$ and $r>1$ then the sequence of partial quotients for $\alpha$ is unbounded. Note that this is the case for continued fractions described by $(S)$. See \cite{S}, for a full study of such elements.

\par A natural question is : what is the size of $\H(q)$ within the algebraic elements of $\F(q)$ ? 
\newline If $\alpha$ is an algebraic element of degree $n$ over $\F_q(T)$ then the formal derivative of $\alpha$ belongs to the vector space $\F_q(T,\alpha)$ of degree $n$ over $\F_q(T)$ and we have
$$\alpha'=a_0+a_1\alpha+a_2\alpha^2+\cdots+a_{n-1}\alpha^{n-1}\quad \text{ where }\quad a_i\in\F_q(T).$$
If $\alpha \in \H(q)$, it is known that $a_i=0$ for $i>2$ ($\alpha$ is differential-quadratic). Hence we see that the probability, for an algebraic element of large degree, to be hyperquadratic is small. For an element of algebraic degree 4, if $a_3\neq 0$ then it is not an hyperquadratic element.
\par A remarkable fact about this subset of algebraic elements is that it contains very well approximable elements as well as badly approximable ones: we can have $\nu(\alpha)=d$ (as for Mahler's example or its dual) or $\nu(\alpha)=2$ (with bounded partial quotients). The continued fraction for Mahler's example is explicitly known \cite[p.~215]{L}, as well as for examples with partial quotients all of degree one. However, the description of the continued fraction expansion for all elements in $\H(q)$ is out of reach.
\par Concerning the algebraic elements which are not hyperquadratic not so much is known. The irrationality measure has been computed for some particular elements, but more rarely the continued fraction (see Section 7, example 1). However in a joint work with de Mathan \cite{LdM} we could prove the following :
\newline Let $\alpha \in \F(q)$ algebraic over $\F_q(T)$ of degree $d>1$. If $\alpha$ is not hyperquadratic then we have $\nu(\alpha)\leq [d/2]+1$.
\newline This result is optimal: if $\alpha$ is Mahler's example, $\alpha$ has degree $r$ and $\nu(\alpha)=r$, then one can prove that $\alpha^2$ has also degree $r$, is not hyperquadratic and we have $\nu(\alpha^2)=r/2$.

\vskip 0.5 cm
\noindent{\bf{7. Three examples of particular algebraic continued fractions}} 
\par We give here three examples of algebraic power series with their full continued fraction expansion.
\par The first two examples come from an algebraic equation of degree 4 introduced by Mills and Robbins \cite[p.~403-404]{MR}: $x^4+x^2-Tx+1=0\quad (MR)$. This equation has a unique root $\alpha$ in $\F(p)$ for all $p$. Note that in all cases we can write :
$$x=(1/T)(1+x^2+x^4)=1/T+u=1/T+1/T((1/T+u)^2+(1/T+u)^4)=\cdots $$
This leads to $\alpha=1/T+1/T^3+\cdots $.
\newline It was observed that the continued fraction for the solution is particularly interesting for $p=3$ and $p=13$.
\newline The last example concerns the solution in $\F(p)$ of an hyperquadratic equation of degree $p+1$, for all $p\geq 3$.
\par {\bf{Example 1:}} Let $(W_n)_{n\geq 0}$ be the sequence of finite words, defined recursively as follows:
$$W_{0}=\epsilon ,\quad W_{1}=T,\quad \textrm{ and }\quad$$ 
$$W_{n}=W_{n-1},2T,W_{n-2}^{(3)},2T,W_{n-1},\quad \text{ for }\quad n\geq 2,$$
where commas indicate concatenation of words, $\epsilon$ denotes the empty word and $W^{(3)}$ denotes the word obtained by rising each letter of $W$ to the power 3. Since $W_n$ begins with  $W_{n-1}$, these sequences approach an infinite word $W_{\infty}(T): T,2T,2T,T,T^3,2T,T,2T,\cdots$.
\newline The continued fraction expansion for the solution $\alpha$ in $\F(3)$ of Mills-Robbins quartic equation is $\alpha=[0,W_{\infty}(T)]$. 
\newline We observe that all the partial quotients are monomials in $\F_3[T]$ of degree a power of three. Moreover the growth of these degrees is very slow, but they are unbounded. A basic computation shows that $$\limsup_{n\geq 1}\deg(a_{n+1})/\sum_{1\leq i\leq n}\deg(a_i)=0.$$ 
\newline {\bf{Exercise 9:}} Prove the above equality.
\newline {\bf{Exercise 10:}} Using the algebraic equation satisfied by $\alpha$, show that it is not differential-quadratic and we have $\alpha'=(\alpha^3-\alpha)/T$.
\newline Using the formula (11) for the irrationality measure, we obtain $\nu(\alpha)=2$.
\newline This fact has an important consequence. Voloch \cite{V} proved that if $\alpha \in \F(q)$ is hyperquadratic then the irrationality measure is 2 if and only if the sequence of the degrees of the partial quotient is bounded. Consequently the solution of $(MR)$ in $\F(3)$ is not hyperquadratic, which is also a direct consequence of Exercise 10. 
\par There are two proofs to obtain this continued fraction expansion (see \cite{BR,L2}). The second one shows a connection with another hyperquadratic element and this may explain why this expansion is not as chaotic as it is for the solution of $(MR)$ with other values of $p$.
\newline Set $\beta=[T,T^3\cdots,T^{3^n},\cdots]\in \F(3)$, note that $\beta$ is the simple hyperquadratic element considered in the previous section. We have $\beta=T+1/\beta^3$. By elevating to the power 2, we get $\beta^2=T^2-T/\beta^3+1/\beta^6$. This becomes $\beta^8-T^2\beta^6+T\beta^3-1=0$. Since $T\beta^3=\beta^4-1$, we get $\beta^8-T^2\beta^6+\beta^4+1=0$ or $(1/\beta^2)^4+(1/\beta^2)^2-T^2(1/\beta^2)+1=0$. this implies $\beta^2=1/\alpha(T^2)$. Consequently we have : 
$$[T,T^3\cdots,T^{3^n},\cdots]^2=[W_{\infty}(T^2)] \quad \text{ in }\quad \F(3).$$
Note that the starting point of this second proof for the continued fraction of the solution of $(MR)$ with $p=3$, was a generalization of the above formula. We wanted to study the following expression in $\F(p)$
$$[T,T^r\cdots,T^{r^n},\cdots]^{(r+1)/2} \quad \text{ where }\quad r=p^t \text{ with } t\geq 1.$$
It is not known, if a continued fraction expansion can be given for this expression for $r>3$ (even in the simplest case $r=5$).
\par {\bf{Example 2:}} The second example is the continued fraction expansion of the solution of $(MR)$ if $p=13$. Here again there is a specificity due to the particular characteristic.
\newline {\bf{Exercise 11:}} Starting from the algebraic equation satisfied by $\alpha$, show that we have $\alpha'=(2\alpha^2-T\alpha/4+1)/(9T^2-6)$. 
\newline Besides, a straight computation shows that $\alpha \in \H(13)$ and we have :
$$9T\alpha^{14}-(T^2+1)\alpha^{13}+(T^6+T^4+11T^2+1)\alpha-(T^5+2T^3+2T)=0.$$
\par Setting $\beta=1/\alpha$, one can prove that $\beta$ is defined by:
$$\beta=[T,12T,7T,11T,8T,5T,\beta_7]\quad \text{ and }\quad $$
$$\beta^{13}=(T^2+8)^4\beta_7+8T^7+T^5+9T^3+7T. \eqno{(12)}$$
In \cite{L3}, we introduced a method to study continued fractions $\alpha$ in $\F(q)$ such that
$$\alpha=[a_1,a_2,\cdots,a_l,\alpha_{l+1}]\quad \text{ and }\quad \alpha^r=P\alpha_{l+1}+Q,$$
where $r=p^t$ and $(P,Q)\in \F_q[T]^2$ is a particular pair. These continued fractions are all hyperquadratic. Indeed, using formula (9) of section 4, we see that $\alpha$ satisfies the following algebraic equation of degree $r+1$ :
$$\alpha=\frac{x_l(\alpha^r-Q)+x_{l-1}P}{y_l(\alpha^r-Q)+y_{l-1}P}. $$
 With a particular choice of $(P,Q)$, generalizing the polynomials appearing in $(12)$, this method was developed in \cite{L4} and the explicit continued fraction for the solution of $(MR)$ in $\F(13)$ was given. The partial quotients $a_n$ are of the form $\lambda_n A_{i(n)}$ where $\lambda_n\in \F_{13}^*$, $i(n)\in \N$ and $(A_k)_{k\geq 0}$ is a particular sequence of polynomials in $\F_{13}[T]$. We indicate here the definition of the sequence $(A_k)_{k\geq 0}$, we have
$$A_0=T\quad \text{ and }\quad A_{k+1}=[A_k^{13}/(T^2+8)^4] \quad \text{ for }\quad k\geq 0.$$
We have $i(n)=v_9(4n-1)$ for $n\geq 1$, where $v_9(m)$ is the largest power of $9$ dividing $m$. The sequence $(\lambda_n)_{n\geq 1}$ in $\F_{13}^*$ is particularly sophisticated and will not be presented here. Note that the irrationality measure is depending on both sequences $(A_k)_{k\geq 0}$ and $(i(n))_{n\geq 1}$, but it is independent of the sequence $(\lambda_n)_{n\geq 1}$. Applying the previous formula, we get here $\nu(\alpha)=8/3 (\in[2;4])$.
\par This second example led to a wide generalization. Later, in\cite{L5}, we extended this study to many more hyperquadratic continued fractions in all odd characteristics, taking $P=(T^2+a)^k$ and $Q$ well chosen. One of the interesting sides of this study is the complexity of the sequence  $(\lambda_n)_{n\geq 1}$ in $\F_q^*$. 

\par {\bf{Example 3:}} This last example is a simple case of a large class of elements in $\F(p)$. These elements are close to those obtained by generalizing the previous example. We consider a prime $p\geq 3$ and a pair $(a,b)\in (\F_p^*)^2$. We shall define an infinite continued fraction $\theta_{a,b}\in \F(p)$, by describing the sequence of partial quotients in $\F_p[T]$.  

\par First, we introduce a sequence $(B_n)_{n\geq 0}$ of unitary polynomials in $\F_p[T]$, defined recursively as follows:
$$B_0=T,\quad B_1=(T^p-T)/(T^2-1)=T\sum_{i=0}^{(p-3)/2}T^{2i} $$
 and
$$B_n=B_{n-1}^p(T^2-1)^{(-1)^n} \quad \text{ for}\quad n\geq 2.$$
\par For $n\in \N^*$, there exists a unique pair of integers $(m,i)$ depending on $n$, such that
$$m\in \N^*, \quad 1\leq i\leq 2m \quad \text{ and }\quad n=m^2-m+i.$$ 
Then, we define a sequence $(\lambda_n)_{n\geq 1}$ in $\F_p^*$ in the following way. For $m\geq 1$ and $1\leq i\leq 2m$, we have 
$$\lambda_n=a  \quad \text{ if }\quad m\geq 1 \quad \text{ and}\quad i=1$$
$$\lambda_n=-b^{(-1)^i}  \quad \text{ if }\quad m\geq 1 \quad \text{ and}\quad  2\leq i\leq m$$
 $$\lambda_n=b^{(-1)^i}  \quad \text{ if }\quad m\geq 1 \quad \text{ and}\quad  m+1\leq i\leq 2m.$$
Hence, we have: $(\lambda_n)_{n\geq 1}=a,b,a,-b,b^{-1},b,a,-b,-b^{-1},b,b^{-1},b,a,\cdots $.
\par Finally, let us consider in $\N$ the sequence $(j(n))_{n\geq 1}$, defined  as follows. For $m\geq 1$ and $1\leq i\leq 2m$, we have 
$$j(n)=m-i \quad \text{ if}\quad 1\leq i\leq m$$
and
$$j(n)=i-m-1 \quad \text{ if}\quad m+1 \leq i\leq 2m.$$
Hence, we have: $(j(n))_{n\geq 1}=0,0,1,0,0,1,2,1,0,0,1,2,3,\cdots $.
\par Then we set
$$\theta_{a,b}=[a_1,a_2,\cdots,a_n,\cdots] \quad \text{where}\quad a_n=\lambda_nB_{j(n)} \quad \text{for}\quad n\geq 1$$
and $\theta_{a,b}$ satisfies the hyperquadratic equation :
$$bTX^{p+1}-(abT^2+1)(X^p+X)+a^2bT^3+(2a+b^{-1})T=0.  \eqno{(PB)}$$ 
It is necessary to underline that the statement just presented above is an \emph{unpublished conjecture !!}
\newline However the rightness is shown by computer calculations. We present here below a short program written with SAGE. This program returns the first partial quotient of the solution of certain algebraic equations with coefficients in $\F_p[T]$. The reader will see here the result when considering $(PB)$ with $p=7$ and $(a,b)=(1,2)$ and he can check that the first partial quotients are conform to the formulas described above.
\newline The connection with example 2 is due to the following: $\alpha\in \F(p)$ satisfies $(PB)$ if and only if we have a definition similar to (12):
$$\alpha=[aT,bT,\alpha_3]\quad \text{and}\quad \alpha^p=(T^2-1)\alpha_3+(a+b^{-1})T.$$
Finally, we can also compute the irrationality measure for $\theta_{a,b}$. Again this irrationality measure is depending on the sequences $(B_n)_{n\geq 0}$ and $(j(n))_{n\geq 1}$, but it is independent of $(\lambda_n)_{n\geq 1}$. We have
$$\nu(\theta_{a,b})=2+(p-1)^2/(2p) \quad \text{for all}\quad p\geq 3.$$
Setting $b_n=\deg(B_n)$, we have $b_n=pb_{n-1}+2(-1)^n$. One can see using formula (11) that we have $\nu(\theta_{a,b})=2+(1/2)\lim_{n\geq 1}b_n/(b_{n-1}+2b_{n-2}+\cdots+nb_0)$. Therefrom we get the value given above.
\newline The exact algebraic degree $d$ of $\theta_{a,b}$ is unknown. Using the irrationality measure and the degree of the equation $(PB)$ satisfied by $\theta_{a,b}$, we get $d\in [(p+3)/2;p+1]$.  
\vskip 0.5 cm
  \centerline{PROGRAMMING WITH SAGE}

\begin{verbatim}
def contf(P,m):
    n=P.degree()
    a=[]
    for i in range(m):
        an=P[n]
        an1=P[n-1]
        qp=-an1//an
        if P(qp)==0:
            i=m
        else:
            P=P(x+qp)
            P=P.reverse()
            
            a.append(qp)
    if n==1:
        a.append(qp)            
    return a

p = 7
F = GF(p)
a=F(1)
b=F(2)
c=2*a+1/b
Ft.<t> = PolynomialRing(F);
Ftx.<x> = PolynomialRing(Ft);
P = b*t*x^(p+1) - (a*b*t^2+1)*(x^p+x) +a^2*b*t^3+c*t
contf(P,12)


[t,
 2*t,
 t^5 + t^3 + t,
 5*t,
 4*t,
 2*t^5 + 2*t^3 + 2*t,
 t^37 + 6*t^35 + t^23 + 6*t^21 + t^9 + 6*t^7,
 5*t^5 + 5*t^3 + 5*t,
 3*t,
 2*t,
 4*t^5 + 4*t^3 + 4*t,
 2*t^37 + 5*t^35 + 2*t^23 + 5*t^21 + 2*t^9 + 5*t^7]

\end{verbatim}

\vskip 1 cm 
\par Using the program given above with SAGE, we have collected several data concerning the examples discussed in this section and other examples. These data are presented in the following annex. 

\par The different examples treated in this section have a common origin. This is explained in a recent note left on Arxiv, we give here the address for the interested reader : https://arxiv.org/abs/1704.08959
\vskip 1 cm
\centerline{\bf{Aknowledgements}}
 \noindent \emph{ These notes were prepared for some classes given at Beihang university (Beijing, P. R. of China) during the summer 2017. I wish to express my gratitude to Professor Zhiyong Zheng for his kind invitation. Let me also warmly thank his young assistant and student, Ziwei Hong, for his valuable and very friendly help during our stay in Beijing. }

\newpage
\centerline{\bf{Annex}}
\centerline{\bf{COMPUTER SCREEN VIEWS }} 

\section{Robbins quartic equation}

$$x^4+x^2-Tx+1=0.$$

We have a solution in $\F(p)$ for all primes $p$. This solution is expanded in continued fraction and we have
$$\alpha=[0,a_1,a_2,a_3,\cdots,a_n,\cdots]$$
Below is the list of the first partial quotients. First, is the list of 10 full partial quotients. Then below is a list of the first 300 partial quotients only given by the degrees.

\vskip 1 cm 
\par $\bullet$ : $p=2$
\begin{verbatim}

 [t, t, t^5 + t, t, t, t, t, t, t^5 + t, t]

degrees [1, 1, 5, 1, 1, 1, 1, 1, 5, 1, 9, 1, 9, 1, 1, 1, 1, 1, 1, 1,

17, 1, 9, 1, 5, 1, 1, 1, 13, 1, 1, 1, 5, 1, 1, 1, 1, 1, 1, 1, 13, 1,

17, 1, 1, 1, 5, 1, 13, 1, 5, 1, 1, 1, 5, 1, 5, 1, 5, 1, 1, 1, 1, 1,

13, 1, 5, 1, 5, 1, 1, 1, 1, 1, 1, 1, 1, 1, 1, 1, 1, 1, 1, 1, 9, 1,

1, 1, 1, 1, 1, 1, 1, 1, 1, 1, 1, 1, 9, 1, 1, 1, 5, 1, 1, 1, 5, 1, 1,

1, 9, 1, 1, 1, 1, 1, 5, 1, 9, 1, 5, 1, 1, 1, 41, 1, 1, 1, 5, 1, 9, 1,

1, 1, 1, 1, 1, 1, 1, 1, 9, 1, 1, 1, 9, 1, 1, 1, 1, 1, 1, 1, 1, 1, 9, 

1, 13, 1, 1, 1, 1, 1, 1, 1, 1, 1, 1, 1, 1, 1, 5, 1, 5, 1, 5, 1, 1, 1,

1, 1, 1, 1, 1, 1, 1, 1, 5, 1, 13, 1, 5, 1, 9, 1, 1, 1, 1, 1, 29, 1, 5,

1, 21, 1, 13, 1, 17, 1, 5, 1, 1, 1, 1, 1, 1, 1, 9, 1, 1, 1, 1, 1, 5, 1,

1, 1, 5, 1, 13, 1, 1, 1, 9, 1, 1, 1, 1, 1, 1, 1, 17, 1, 1, 1, 5, 1, 1,

1, 13, 1, 5, 1, 1, 1, 9, 1, 13, 1, 1, 1, 1, 1, 5, 1, 9, 1, 9, 1, 5, 1,

1, 1, 5, 1, 5, 1, 5, 1, 1, 1, 5, 1, 65, 1, 9, 1, 1, 1, 1, 1, 5, 1, 5, 

1, 1, 1, 1, 1, 1, 1]


\end{verbatim}
\vskip 1 cm 
\par $\bullet$ : $p=3$
\begin{verbatim}


[t, 2*t, 2*t, t, 2*t, t^3, 2*t, t, 2*t, 2*t]

degrees [1, 1, 1, 1, 1, 3, 1, 1, 1, 1, 1, 1, 3, 3, 3, 3, 1, 1,

1, 1, 1, 1, 3, 1, 1, 1, 1, 1, 1, 3, 3, 3, 3, 3, 9, 3, 3, 3, 3,

3, 1, 1, 1, 1, 1, 1, 3, 1, 1, 1, 1, 1, 1, 3, 3, 3, 3, 1, 1, 1,

1, 1, 1, 3, 1, 1, 1, 1, 1, 1, 3, 3, 3, 3, 3, 9, 3, 3, 3, 3, 3,

3, 9, 9, 9, 9, 3, 3, 3, 3, 3, 3, 9, 3, 3, 3, 3, 3, 1, 1, 1, 1,

1, 1, 3, 1, 1, 1, 1, 1, 1, 3, 3, 3, 3, 1, 1, 1, 1, 1, 1, 3, 1,

1, 1, 1, 1, 1, 3, 3, 3, 3, 3, 9, 3, 3, 3, 3, 3, 1, 1, 1, 1, 1,

1, 3, 1, 1, 1, 1, 1, 1, 3, 3, 3, 3, 1, 1, 1, 1, 1, 1, 3, 1, 1,

1, 1, 1, 1, 3, 3, 3, 3, 3, 9, 3, 3, 3, 3, 3, 3, 9, 9, 9, 9, 3,

3, 3, 3, 3, 3, 9, 3, 3, 3, 3, 3, 3, 9, 9, 9, 9, 9, 27, 9, 9, 9,

9, 9, 3, 3, 3, 3, 3, 3, 9, 3, 3, 3, 3, 3, 3, 9, 9, 9, 9, 3, 3,

3, 3, 3, 3, 9, 3, 3, 3, 3, 3, 1, 1, 1, 1, 1, 1, 3, 1, 1, 1, 1,

1, 1, 3, 3, 3, 3, 1, 1, 1, 1, 1, 1, 3, 1, 1, 1, 1, 1, 1, 3, 3,

3, 3, 3, 9, 3, 3, 3, 3, 3, 1, 1, 1, 1, 1, 1, 3, 1, 1, 1, 1, 1,

1, 3, 3, 3, 3, 1, 1, 1, 1]


\end{verbatim}
\vskip 1 cm 
\par $\bullet$ : $p=5$
\begin{verbatim}

[t, 4*t, 3*t, 3*t, 4*t^3 + 4*t, 2*t, 3*t^3 + t, 3*t, 3*t^5 + 3*t^3, 4*t^3 + 2*t]

degrees [1, 1, 1, 1, 3, 1, 3, 1, 5, 3, 3, 1, 1, 1, 1, 1, 1, 3, 3, 1, 1,

1, 3, 1, 3, 1, 1, 3, 1, 1, 1, 3, 5, 1, 1, 1, 1, 1, 1, 1, 1, 3, 1, 3, 1,

1, 1, 1, 3, 1, 5, 1, 1, 1, 5, 1, 1, 1, 1, 1, 1, 1, 1, 1, 1, 1, 1, 3, 1,

1, 3, 1, 1, 1, 1, 3, 1, 1, 1, 3, 1, 3, 1, 1, 1, 1, 1, 5, 3, 1, 1, 1, 1,

1, 1, 3, 3, 3, 1, 1, 3, 3, 3, 1, 1, 1, 3, 1, 1, 5, 1, 1, 1, 1, 3, 3, 3,

1, 1, 1, 1, 1, 1, 1, 3, 1, 3, 1, 1, 3, 1, 1, 1, 1, 3, 1, 1, 1, 1, 3, 1,

1, 1, 3, 1, 1, 1, 3, 3, 1, 1, 1, 1, 1, 1, 1, 1, 1, 1, 1, 3, 3, 1, 1, 5,

1, 1, 1, 1, 1, 1, 1, 1, 1, 1, 1, 5, 1, 3, 1, 1, 1, 1, 3, 1, 1, 1, 1, 1,

1, 1, 1, 1, 1, 1, 1, 5, 1, 1, 3, 7, 1, 1, 1, 1, 1, 3, 3, 1, 3, 1, 1, 1,

3, 1, 1, 1, 1, 1, 1, 5, 1, 3, 1, 1, 1, 5, 1, 1, 3, 1, 1, 3, 1, 1, 1, 3, 

1, 1, 1, 1, 3, 1, 1, 1, 1, 1, 1, 1, 1, 1, 1, 3, 1, 1, 1, 3, 1, 1, 1, 3,

1, 1, 1, 1, 1, 3, 1, 1, 1, 1, 1, 1, 3, 1, 1, 1, 1, 3, 1, 1, 1, 3, 1, 1,

1, 3, 1, 1, 7, 1, 1, 1, 1, 3, 1, 1, 3, 1, 1]



\end{verbatim}
\vskip 1 cm 
\par $\bullet$ : $p=7$
\begin{verbatim}

[t, 6*t, 4*t, 5*t, 3*t, 2*t, t, 2*t, t, 2*t]

degrees [1, 1, 1, 1, 1, 1, 1, 1, 1, 1, 1, 1, 1, 1,

1, 1, 1, 1, 1, 3, 1, 1, 3, 1, 1, 1, 1, 1, 1, 1, 1,

1, 1, 1, 1, 1, 1, 3, 3, 1, 1, 1, 1, 1, 1, 1, 1, 1,

1, 1, 1, 3, 1, 1, 1, 3, 1, 1, 1, 1, 1, 1, 1, 3, 1,

1, 1, 5, 1, 1, 1, 1, 1, 1, 1, 3, 1, 1, 1, 1, 1, 1,

1, 1, 1, 1, 1, 3, 1, 1, 3, 3, 3, 1, 1, 7, 1, 3, 1,

3, 1, 1, 1, 1, 1, 3, 1, 1, 1, 1, 1, 3, 1, 1, 1, 1, 

1, 1, 3, 1, 3, 1, 1, 1, 1, 1, 1, 1, 3, 1, 1, 1, 1, 

1, 1, 1, 1, 1, 3, 1, 1, 1, 1, 3, 3, 3, 1, 1, 1, 1,

1, 1, 3, 1, 1, 1, 1, 1, 1, 1, 1, 1, 1, 3, 1, 1, 1,

1, 1, 5, 1, 3, 1, 1, 1, 1, 1, 1, 1, 1, 1, 1, 1, 1,

1, 1, 1, 1, 1, 1, 1, 1, 1, 1, 3, 1, 1, 1, 1, 1, 1,

1, 1, 1, 3, 1, 1, 3, 1, 1, 1, 1, 1, 1, 1, 1, 1, 1,

3, 1, 1, 1, 1, 1, 1, 1, 1, 1, 3, 1, 1, 1, 3, 1, 1,

1, 1, 1, 1, 1, 3, 1, 1, 1, 1, 1, 1, 1, 1, 1, 1, 1,

1, 1, 1, 1, 3, 1, 1, 1, 1, 1, 1, 1, 1, 1, 1, 1, 1,

1, 1, 1, 1, 1, 1, 1, 1, 3, 1, 1, 3, 1, 1, 1, 1, 1,

3, 1, 1, 3, 1, 3, 1, 1, 1, 1, 3, 1, 3, 1]

\end{verbatim}

\vskip 1 cm 
\par $\bullet$ : $p=11$
\begin{verbatim}
 [t, 10*t, 6*t, 9*t, 9*t, 9*t, 3*t^3 + 5*t, 8*t, 2*t, 3*t]
 
degrees [1, 1, 1, 1, 1, 1, 3, 1, 1, 1, 1, 1, 1, 1, 1,

1, 1, 1, 3, 1, 1, 1, 1, 1, 3, 3, 1, 1, 1, 3, 1, 1, 3, 

1, 1, 1, 1, 3, 1, 3, 1, 1, 1, 1, 1, 1, 1, 3, 1, 1, 3,

1, 1, 1, 1, 1, 1, 1, 1, 1, 1, 1, 1, 1, 1, 1, 1, 1, 1,

1, 1, 1, 1, 1, 1, 1, 1, 1, 1, 1, 3, 1, 1, 1, 1, 1, 1,

3, 1, 1, 1, 1, 1, 1, 1, 1, 1, 1, 1, 1, 1, 1, 1, 1, 1,

1, 1, 1, 1, 1, 1, 1, 1, 3, 3, 1, 3, 1, 1, 3, 1, 1, 1,

1, 1, 1, 1, 1, 1, 1, 1, 3, 1, 1, 1, 1, 1, 1, 1, 1, 1,

1, 3, 1, 1, 1, 1, 1, 1, 3, 1, 1, 1, 1, 1, 1, 1, 1, 1,

1, 1, 3, 1, 1, 1, 1, 1, 1, 1, 3, 1, 1, 1, 1, 1, 1, 1,

1, 1, 1, 1, 1, 1, 1, 1, 1, 1, 1, 1, 1, 1, 1, 1, 1, 1, 

1, 1, 1, 1, 1, 1, 1, 1, 5, 1, 1, 1, 1, 1, 1, 1, 1, 1,

1, 1, 3, 3, 1, 1, 3, 1, 1, 1, 1, 1, 1, 3, 1, 1, 1, 1,

1, 1, 1, 1, 1, 1, 1, 1, 1, 1, 3, 1, 1, 1, 1, 1, 1, 1,

1, 1, 3, 1, 1, 1, 1, 1, 1, 1, 1, 1, 3, 1, 1, 1, 1, 1,

1, 1, 1, 1, 1, 1, 1, 3, 1, 1, 1, 1, 1, 1, 3, 1, 1, 1,

1, 1, 1, 1, 1, 1, 1, 1, 1, 1, 3, 1, 1, 1, 1]


\end{verbatim}
\vskip 1 cm 
\par $\bullet$ : $p=13$
\begin{verbatim}

[t, 12*t, 7*t, 11*t, 8*t, 5*t, t^5 + 7*t^3 + 3*t, 3*t, 9*t, 4*t]

degrees [1, 1, 1, 1, 1, 1, 5, 1, 1, 1, 1, 1, 1, 1, 1, 5,

1, 1, 1, 1, 1, 1, 1, 1, 5, 1, 1, 1, 1, 1, 1, 1, 1, 5, 1,

1, 1, 1, 1, 1, 1, 1, 5, 1, 1, 1, 1, 1, 1, 1, 1, 5, 1, 1,

1, 1, 1, 1, 1, 1, 57, 1, 1, 1, 1, 1, 1, 1, 1, 5, 1, 1, 1,

1, 1, 1, 1, 1, 5, 1, 1, 1, 1, 1, 1, 1, 1, 5, 1, 1, 1, 1,

1, 1, 1, 1, 5, 1, 1, 1, 1, 1, 1, 1, 1, 5, 1, 1, 1, 1, 1,

1, 1, 1, 5, 1, 1, 1, 1, 1, 1, 1, 1, 5, 1, 1, 1, 1, 1, 1,

1, 1, 5, 1, 1, 1, 1, 1, 1, 1, 1, 57, 1, 1, 1, 1, 1, 1, 1, 

1, 5, 1, 1, 1, 1, 1, 1, 1, 1, 5, 1, 1, 1, 1, 1, 1, 1, 1, 

5, 1, 1, 1, 1, 1, 1, 1, 1, 5, 1, 1, 1, 1, 1, 1, 1, 1, 5,

1, 1, 1, 1, 1, 1, 1, 1, 5, 1, 1, 1, 1, 1, 1, 1, 1, 5, 1,

1, 1, 1, 1, 1, 1, 1, 5, 1, 1, 1, 1, 1, 1, 1, 1, 57, 1, 1,

1, 1, 1, 1, 1, 1, 5, 1, 1, 1, 1, 1, 1, 1, 1, 5, 1, 1, 1,

1, 1, 1, 1, 1, 5, 1, 1, 1, 1, 1, 1, 1, 1, 5, 1, 1, 1, 1,

1, 1, 1, 1, 5, 1, 1, 1, 1, 1, 1, 1, 1, 5, 1, 1, 1, 1, 1,

1, 1, 1, 5, 1, 1, 1, 1, 1, 1, 1, 1, 5, 1, 1, 1, 1, 1]


\end{verbatim}

\section{Modified Robbins quartic equation}

$$x^4+x^2-Tx-1/12=0\quad \text{ for} \quad p>3.$$

We have a solution in $\F(p)$ for all primes $p>3$. This solution is expanded in continued fraction and we have
$$\alpha=[0,a_1,a_2,a_3,\cdots,a_n,\cdots]$$
It was proved that this solution is hyperquadratic of order $1$ if $p\equiv 1 \mod 3$ and hyperquadratic of order $2$ if $p\equiv 2 \mod 3$.
\newline Here below we show the list of the first 300 partial quotients (only given by their degrees) for the first four values of $p$ : $5,7,11$ and $13$ (also a list of the first partial quotients and secondly of their leading coefficients).
\vskip 0.5 cm 
\par $\bullet$ : $p=5$
\begin{verbatim}
cfe [3*t, 4*t, t, t, t, t, 4*t, 3*t, 2*t, 3*t]
degrees [1, 1, 1, 1, 1, 1, 1, 1, 1, 1, 1, 1, 9, 41, 9, 1,

1, 1, 1, 1, 1, 1, 1, 1, 1, 1, 1, 1, 1, 1, 1, 9, 1, 1, 1,

1, 1, 1, 1, 1, 1, 1, 1, 1, 1, 1, 1, 1, 9, 1, 1, 1, 1, 1, 

1, 1, 1, 1, 1, 1, 1, 1, 1, 1, 1, 9, 41, 9, 1, 1, 1, 1, 1,

1, 1, 1, 1, 1, 1, 1, 1, 1, 1, 1, 9, 1, 1, 1, 1, 1, 1, 1, 

1, 1, 1, 1, 1, 1, 1, 1, 1, 9, 1, 1, 1, 1, 1, 1, 1, 1, 1, 

1, 1, 1, 1, 1, 1, 1, 9, 41, 209, 1041, 209, 41, 9, 1, 1, 

1, 1, 1, 1, 1, 1, 1, 1, 1, 1, 1, 1, 1, 1, 9, 1, 1, 1, 1, 

1, 1, 1, 1, 1, 1, 1, 1, 1, 1, 1, 1, 9, 1, 1, 1, 1, 1, 1,

1, 1, 1, 1, 1, 1, 1, 1, 1, 1, 9, 41, 9, 1, 1, 1, 1, 1, 1,

1, 1, 1, 1, 1, 1, 1, 1, 1, 1, 9, 1, 1, 1, 1, 1, 1, 1, 1, 

1, 1, 1, 1, 1, 1, 1, 1, 9, 1, 1, 1, 1, 1, 1, 1, 1, 1, 1, 

1, 1, 1, 1, 1, 1, 9, 41, 9, 1, 1, 1, 1, 1, 1, 1, 1, 1, 1,

1, 1, 1, 1, 1, 1, 9, 1, 1, 1, 1, 1, 1, 1, 1, 1, 1, 1, 1,

1, 1, 1, 1, 9, 1, 1, 1, 1, 1, 1, 1, 1, 1, 1, 1, 1, 1, 1, 

1, 1, 9, 41, 209, 41, 9, 1, 1, 1, 1, 1, 1, 1, 1, 1, 1, 1,

1, 1, 1]

leading coefficients
[3, 4, 1, 1, 1, 1, 4, 3, 2, 3, 2, 4, 2, 1, 4, 2, 4, 4, 4, 4]
\end{verbatim}
\vskip 1 cm 
\par $\bullet$ : $p=7$
\begin{verbatim}
cfe [2*t, 6*t, 6*t, 3*t^3 + 6*t, 5*t, 3*t, 4*t, 2*t, 4*t^3 + t, t]

degrees [1, 1, 1, 3, 1, 1, 1, 1, 3, 1, 1, 1, 1, 3, 1, 1, 1, 1, 17,

1, 1, 1, 1, 3, 1, 1, 1, 1, 3, 1, 1, 1, 1, 3, 1, 1, 1, 1, 3, 1, 1,

1, 1, 17, 1, 1, 1, 1, 3, 1, 1, 1, 1, 3, 1, 1, 1, 1, 3, 1, 1, 1, 1,

3, 1, 1, 1, 1, 17, 1, 1, 1, 1, 3, 1, 1, 1, 1, 3, 1, 1, 1, 1, 3, 1,

1, 1, 1, 3, 1, 1, 1, 1, 115, 1, 1, 1, 1, 3, 1, 1, 1, 1, 3, 1, 1, 1,

1, 3, 1, 1, 1, 1, 3, 1, 1, 1, 1, 17, 1, 1, 1, 1, 3, 1, 1, 1, 1, 3,

1, 1, 1, 1, 3, 1, 1, 1, 1, 3, 1, 1, 1, 1, 17, 1, 1, 1, 1, 3, 1, 1,

1, 1, 3, 1, 1, 1, 1, 3, 1, 1, 1, 1, 3, 1, 1, 1, 1, 17, 1, 1, 1, 1,

3, 1, 1, 1, 1, 3, 1, 1, 1, 1, 3, 1, 1, 1, 1, 3, 1, 1, 1, 1, 17, 1,

1, 1, 1, 3, 1, 1, 1, 1, 3, 1, 1, 1, 1, 3, 1, 1, 1, 1, 3, 1, 1, 1,

1, 115, 1, 1, 1, 1, 3, 1, 1, 1, 1, 3, 1, 1, 1, 1, 3, 1, 1, 1, 1, 

3, 1, 1, 1, 1, 17, 1, 1, 1, 1, 3, 1, 1, 1, 1, 3, 1, 1, 1, 1, 3, 1,

1, 1, 1, 3, 1, 1, 1, 1, 17, 1, 1, 1, 1, 3, 1, 1, 1, 1, 3, 1, 1, 1,

1, 3, 1, 1, 1, 1, 3, 1, 1, 1, 1, 17, 1, 1, 1, 1, 3, 1]

leading coefficients
[2, 6, 6, 3, 5, 3, 4, 2, 4, 1, 1, 5, 3, 2, 1, 1, 5, 3, 2, 6]
\end{verbatim}

\noindent {\bf{Exercise 12:}} Let $\beta=1/\alpha$. Show that we have
$$(2T^2+2)\beta^8+(3T^3+5T)\beta^7+4T\beta+1=0.$$
Using this last equation and also $\beta=[2T,6T,6T,\beta_4]$, prove that we have
$$\beta^7=3(T^2-1)^2\beta_4+4T^3+2T.$$
Finally, using Exercise 3, give $a_4$ to $a_8$ and prove that we also have
$$\beta_2^7=-(T^2-1)^2\beta_9+3T^3+5T.$$
\vskip 1 cm 
\par $\bullet$ : $p=11$
\begin{verbatim}
cfe [10*t, 10*t, 5*t^3 + t, 7*t, 8*t, 10*t, 7*t, t, 4*t, 9*t]

degrees [1, 1, 3, 1, 1, 1, 1, 1, 1, 1, 1, 3, 1, 1, 1, 1, 1, 1, 1,

1, 3, 1, 1, 1, 1, 1, 1, 1, 1, 3, 1, 1, 1, 1, 1, 1, 1, 1, 3, 1, 1,

1, 1, 1, 1, 1, 1, 3, 41, 3, 1, 1, 1, 1, 1, 1, 1, 1, 3, 1, 1, 1, 

1, 1, 1, 1, 1, 3, 1, 1, 1, 1, 1, 1, 1, 1, 3, 1, 1, 1, 1, 1, 1, 1,

1, 3, 1, 1, 1, 1, 1, 1, 1, 1, 3, 1, 1, 1, 1, 1, 1, 1, 1, 3, 1, 1,

1, 1, 1, 1, 1, 1, 3, 41, 443, 41, 3, 1, 1, 1, 1, 1, 1, 1, 1, 3, 

1, 1, 1, 1, 1, 1, 1, 1, 3, 1, 1, 1, 1, 1, 1, 1, 1, 3, 1, 1, 1, 1,

1, 1, 1, 1, 3, 1, 1, 1, 1, 1, 1, 1, 1, 3, 1, 1, 1, 1, 1, 1, 1, 1,

3, 1, 1, 1, 1, 1, 1, 1, 1, 3, 41, 3, 1, 1, 1, 1, 1, 1, 1, 1, 3, 1,

1, 1, 1, 1, 1, 1, 1, 3, 1, 1, 1, 1, 1, 1, 1, 1, 3, 1, 1, 1, 1, 1,

1, 1, 1, 3, 1, 1, 1, 1, 1, 1, 1, 1, 3, 1, 1, 1, 1, 1, 1, 1, 1, 3,

1, 1, 1, 1, 1, 1, 1, 1, 3, 41, 3, 1, 1, 1, 1, 1, 1, 1, 1, 3, 1, 1,

1, 1, 1, 1, 1, 1, 3, 1, 1, 1, 1, 1, 1, 1, 1, 3, 1, 1, 1, 1, 1, 1,

1, 1, 3, 1, 1, 1, 1, 1, 1, 1, 1, 3, 1, 1, 1, 1, 1, 1, 1, 1]

leading coefficients
[10, 10, 5, 7, 8, 10, 7, 1, 4, 9, 5, 1, 10, 10, 8, 6, 3, 5, 5, 9]

\end{verbatim}
\vskip 1 cm 
\par $\bullet$ : $p=13$
\begin{verbatim}
cfe [t, 12*t, 7*t, 11*t, 8*t, 5*t, t^5 + 7*t^3 + 3*t, 3*t, 9*t, 4*t]

degrees [1, 1, 1, 1, 1, 1, 5, 1, 1, 1, 1, 1, 1, 1, 1, 5, 1, 1, 1, 1,

1, 1, 1, 1, 5, 1, 1, 1, 1, 1, 1, 1, 1, 5, 1, 1, 1, 1, 1, 1, 1, 1, 

5, 1, 1, 1, 1, 1, 1, 1, 1, 5, 1, 1, 1, 1, 1, 1, 1, 1, 57, 1, 1, 1,

1, 1, 1, 1, 1, 5, 1, 1, 1, 1, 1, 1, 1, 1, 5, 1, 1, 1, 1, 1, 1, 1, 

1, 5, 1, 1, 1, 1, 1, 1, 1, 1, 5, 1, 1, 1, 1, 1, 1, 1, 1, 5, 1, 1,

1, 1, 1, 1, 1, 1, 5, 1, 1, 1, 1, 1, 1, 1, 1, 5, 1, 1, 1, 1, 1, 1,

1, 1, 5, 1, 1, 1, 1, 1, 1, 1, 1, 57, 1, 1, 1, 1, 1, 1, 1, 1, 5, 1,

1, 1, 1, 1, 1, 1, 1, 5, 1, 1, 1, 1, 1, 1, 1, 1, 5, 1, 1, 1, 1, 1,

1, 1, 1, 5, 1, 1, 1, 1, 1, 1, 1, 1, 5, 1, 1, 1, 1, 1, 1, 1, 1, 5,

1, 1, 1, 1, 1, 1, 1, 1, 5, 1, 1, 1, 1, 1, 1, 1, 1, 5, 1, 1, 1, 1,

1, 1, 1, 1, 57, 1, 1, 1, 1, 1, 1, 1, 1, 5, 1, 1, 1, 1, 1, 1, 1, 1,

5, 1, 1, 1, 1, 1, 1, 1, 1, 5, 1, 1, 1, 1, 1, 1, 1, 1, 5, 1, 1, 1, 

1, 1, 1, 1, 1, 5, 1, 1, 1, 1, 1, 1, 1, 1, 5, 1, 1, 1, 1, 1, 1, 1,

1, 5, 1, 1, 1, 1, 1, 1, 1, 1, 5, 1, 1, 1, 1, 1]

leading coefficients
[1, 12, 7, 11, 8, 5, 1, 3, 9, 4, 3, 2, 6, 6, 11, 12, 9, 3, 12, 1]
\end{verbatim}

There are clearly two different patterns : cases $p=5$ and $p=11$ (more generally $p\equiv 2 \mod 3$) on one side and  cases $p=7$ and $p=13$ (more generally $p\equiv 1 \mod 3$) on the other side.
\newline In each case ($p\equiv 1 \mod 3$ or $p\equiv 2 \mod 3$) the continued fraction is based on a particular polynomial in $\F_p[T]$ : $P_{u,k}=(T^2+u)^k$ where $u\in \F_p^*$ and $k\in \N$ are two parameters depending on the prime $p$. In the first case, $p\equiv 1 \mod 3$, this means that the continued fraction $\beta=1/\alpha$ is such that we have
$$\beta=[a_1,a_2,\cdots,a_\ell,\beta_{\ell+1}] \quad \text{and}\quad \beta^p=\epsilon_1P_{u,k}\beta_{\ell+1}+\epsilon_2R$$
where $l\geq 1$ is an integer, the pair $(\epsilon_1,\epsilon_2)\in (\F_p*)^2$, all depending on $p$, and $R\in \F_p[T]$ is the remainder in the Euclidean division of $T^p$ by $P_{u,k}$. Besides, the $\ell$ first partial quotients are given and , in this case $p\equiv 1 \mod 3$, these are of the form $a_i=\lambda_iT$ with $\lambda_i\in \F_p^*$. In both cases, all the partial quotients, up to a constant factor, belong to a particular sequence (different according to the case considered) only depending on the polynomial $P_{u,k}$.
\newline The case $p\equiv 1 \mod 3$ enters in a much larger family which has been fully studied. While the second case enters into another family which has only been partially described (see \cite{AL}).
\section{Generalization of case $p=3$ in Robbins quartic equation}

As usual we have $r=p^t$ with $t\geq 1$ and $p$ odd.
\newline We set $$\alpha=[T,T^r,T^{r^2},\cdots,T^{r^n},\cdots]\quad \text{and}\quad \beta=\alpha^{(r+1)/2}.$$
We have $\alpha=T+1/\alpha^r$ and therefrom $T\alpha^r=\alpha^{r+1}-1$. Since $\beta^2=\alpha^{r+1}$, by elevating this formula to the power $(r+1)/2$, we get
$$(\beta^2-1)^{(r+1)/2}-T^{(r+1)/2}\beta^r=0.$$
This formula shows that $\beta$ is an algebraic power series in the variable $T^{(r+1)/2}$. We set $\beta(T)=\gamma(T^{(r+1)/2})$, and $\gamma$ satisfies
$$(\gamma^2-1)^{(r+1)/2}-T\gamma^r=0.$$
Note that for $r=3$ we obtain :$\gamma^4-T\gamma^3 +\gamma^2+1=0$ Hence $1/gamma$ is the solution of robbins quartic equation for $p=3$.
\newline To observe a possible analogy between the different cases $r=3$ and $r=5$, we have given below the list of the first 300 partial quotients for the continued fraction of $\gamma)$.

\vskip 1 cm 
\par $\bullet$ : $r=3$
\begin{verbatim}
cfe [t, 2*t, 2*t, t, 2*t, t^3, 2*t, t, 2*t, 2*t]

degrees [1, 1, 1, 1, 1, 3, 1, 1, 1, 1, 1, 1, 3, 3, 3, 3, 1, 1,

1, 1, 1, 1, 3, 1, 1, 1, 1, 1, 1, 3, 3, 3, 3, 3, 9, 3, 3, 3, 3,

3, 1, 1, 1, 1, 1, 1, 3, 1, 1, 1, 1, 1, 1, 3, 3, 3, 3, 1, 1, 1,

1, 1, 1, 3, 1, 1, 1, 1, 1, 1, 3, 3, 3, 3, 3, 9, 3, 3, 3, 3, 3,

3, 9, 9, 9, 9, 3, 3, 3, 3, 3, 3, 9, 3, 3, 3, 3, 3, 1, 1, 1, 1,

1, 1, 3, 1, 1, 1, 1, 1, 1, 3, 3, 3, 3, 1, 1, 1, 1, 1, 1, 3, 1,

1, 1, 1, 1, 1, 3, 3, 3, 3, 3, 9, 3, 3, 3, 3, 3, 1, 1, 1, 1, 1,

1, 3, 1, 1, 1, 1, 1, 1, 3, 3, 3, 3, 1, 1, 1, 1, 1, 1, 3, 1, 1,

1, 1, 1, 1, 3, 3, 3, 3, 3, 9, 3, 3, 3, 3, 3, 3, 9, 9, 9, 9, 3, 

3, 3, 3, 3, 3, 9, 3, 3, 3, 3, 3, 3, 9, 9, 9, 9, 9, 27, 9, 9, 9,

9, 9, 3, 3, 3, 3, 3, 3, 9, 3, 3, 3, 3, 3, 3, 9, 9, 9, 9, 3, 3, 

3, 3, 3, 3, 9, 3, 3, 3, 3, 3, 1, 1, 1, 1, 1, 1, 3, 1, 1, 1, 1,

1, 1, 3, 3, 3, 3, 1, 1, 1, 1, 1, 1, 3, 1, 1, 1, 1, 1, 1, 3, 3,

3, 3, 3, 9, 3, 3, 3, 3, 3, 1, 1, 1, 1, 1, 1, 3, 1, 1, 1, 1, 1,

1, 3, 3, 3, 3, 1, 1, 1, 1]

leading coefficients 
[1, 2, 2, 1, 2, 1, 2, 1, 2, 2, 1, 2, 1, 2, 2, 1, 2, 1, 2, 2]

\end{verbatim}
\vskip 1 cm 
\par $\bullet$ : $r=5$
\begin{verbatim}
cfe [t, 2*t, 2*t, 2*t, 2*t, t, 2*t, 4*t^3 + 4*t, 4*t, t]

degrees [1, 1, 1, 1, 1, 1, 1, 3, 1, 1, 1, 5, 1, 1, 1, 3, 1, 1,

1, 1, 1, 1, 1, 1, 3, 1, 1, 3, 1, 1, 1, 1, 1, 1, 3, 1, 1, 3, 1,

1, 1, 3, 1, 1, 1, 5, 5, 5, 5, 5, 5, 1, 1, 1, 3, 1, 1, 1, 3, 1,

1, 3, 1, 1, 1, 1, 1, 1, 3, 1, 1, 3, 1, 1, 1, 1, 1, 1, 1, 1, 3,

1, 1, 1, 5, 1, 1, 1, 3, 1, 1, 1, 1, 1, 1, 1, 1, 3, 1, 1, 3, 1,

1, 1, 1, 1, 1, 3, 1, 1, 3, 1, 1, 1, 3, 1, 1, 3, 1, 1, 1, 1, 1, 

1, 1, 5, 1, 1, 3, 1, 1, 5, 1, 1, 3, 1, 1, 3, 1, 1, 1, 1, 3, 1,

3, 1, 3, 1, 1, 1, 1, 1, 1, 3, 1, 7, 1, 3, 1, 1, 1, 3, 1, 1, 1,

1, 1, 1, 1, 3, 3, 1, 1, 1, 3, 1, 1, 3, 1, 1, 1, 1, 1, 1, 3, 1,

1, 3, 1, 1, 1, 3, 1, 1, 1, 5, 5, 5, 5, 5, 5, 5, 15, 5, 5, 5, 25,

5, 5, 5, 15, 5, 5, 5, 5, 5, 5, 5, 1, 1, 1, 3, 1, 1, 1, 3, 1, 1,

3, 1, 1, 1, 1, 1, 1, 3, 1, 1, 3, 1, 1, 1, 3, 3, 1, 1, 1, 1, 1, 1,

1, 3, 1, 1, 1, 3, 1, 7, 1, 3, 1, 1, 1, 1, 1, 1, 3, 1, 3, 1, 3, 1,

1, 1, 1, 3, 1, 1, 3, 1, 1, 5, 1, 1, 3, 1, 1, 5, 1, 1, 1, 1, 1, 1,

1, 3, 1, 1, 3, 1]


leading coefficients
[1, 2, 2, 2, 2, 1, 2, 4, 4, 1, 2, 1, 2, 1, 4, 4, 2, 1, 2, 2]

\end{verbatim}

\section{A last example of a perfect expansion in $\F(5)$}

Here we just give the list as previously for a particular $\alpha$ defined by $$\alpha=[T,T,T,\alpha_4]\quad \text{ and }\quad \alpha^5=(T^2-1)^2\alpha_4+3T^3+T.$$

\vskip 1 cm
\begin{verbatim}
cfe [t, t, t, t, 4*t, 4*t, 4*t, 4*t, t, 2*t]

degrees [1, 1, 1, 1, 1, 1, 1, 1, 1, 1, 1, 1, 1, 1, 9, 1, 1, 1, 1, 

1, 1, 1, 1, 1, 1, 1, 1, 1, 1, 1, 1, 9, 1, 1, 1, 1, 1, 1, 1, 1, 1,

1, 1, 1, 1, 1, 1, 1, 9, 41, 9, 1, 1, 1, 1, 1, 1, 1, 1, 1, 1, 1, 1,

1, 1, 1, 1, 9, 1, 1, 1, 1, 1, 1, 1, 1, 1, 1, 1, 1, 1, 1, 1, 1, 9, 

1, 1, 1, 1, 1, 1, 1, 1, 1, 1, 1, 1, 1, 1, 1, 1, 9, 41, 9, 1, 1, 1,

1, 1, 1, 1, 1, 1, 1, 1, 1, 1, 1, 1, 1, 9, 1, 1, 1, 1, 1, 1, 1, 1,

1, 1, 1, 1, 1, 1, 1, 1, 9, 1, 1, 1, 1, 1, 1, 1, 1, 1, 1, 1, 1, 1,

1, 1, 1, 9, 41, 209, 41, 9, 1, 1, 1, 1, 1, 1, 1, 1, 1, 1, 1, 1, 1,

1, 1, 1, 9, 1, 1, 1, 1, 1, 1, 1, 1, 1, 1, 1, 1, 1, 1, 1, 1, 9, 1, 1,

1, 1, 1, 1, 1, 1, 1, 1, 1, 1, 1, 1, 1, 1, 9, 41, 9, 1, 1, 1, 1, 1, 1,

1, 1, 1, 1, 1, 1, 1, 1, 1, 1, 9, 1, 1, 1, 1, 1, 1, 1, 1, 1, 1, 1, 1,

1, 1, 1, 1, 9, 1, 1, 1, 1, 1, 1, 1, 1, 1, 1, 1, 1, 1, 1, 1, 1, 9, 41,

9, 1, 1, 1, 1, 1, 1, 1, 1, 1, 1, 1, 1, 1, 1, 1, 1, 9, 1, 1, 1, 1, 1,

1, 1, 1, 1, 1, 1, 1, 1, 1, 1, 1, 9, 1]


leading coefficients 
 [1, 1, 1, 1, 4, 4, 4, 4, 1, 2, 3, 2, 3, 1, 1, 4, 2, 3, 2,

3, 4, 1, 1, 1, 1, 4, 3, 2, 3, 2, 4, 1, 2, 4, 4, 4, 4, 3, 3, 2, 3, 2, 

2, 1, 1, 1, 1, 3, 1, 1, 4, 2, 4, 4, 4, 4, 3, 3, 2, 3, 2, 2, 1, 1, 1,

1, 3, 4, 1, 3, 2, 3, 2, 1, 4, 4, 4, 4, 1, 2, 3, 2, 3, 1, 4, 3, 1, 1,

1, 1, 2, 2, 3, 2, 3, 3, 4, 4, 4, 4, 2, 4, 1, 2, 4, 2, 3, 2, 3, 4, 1,

1, 1, 1, 4, 3, 2, 3, 2, 4]

\end{verbatim}

This last example belongs to the same large family as the one presented at the end of Section 7. Up to a constant factor in $\F_5^*$, the partial quotients are known to belong to a special sequence of unitary polynomials. The distribution of these polynomials in the list of partial quotients is understood and this implies that, in this particular case, the irrationality measure of $\alpha$ is equal to $18/7$. However the sequence of the leading coefficients, (see the last table of leading coefficients), is yet difficult to describe.

\vskip 0.5 cm
\begin{tabular}{ll}Alain LASJAUNIAS\\Institut de Math\'ematiques de Bordeaux
\\Universit\'e de Bordeaux, France \\E-mail: Alain.Lasjaunias@math.u-bordeaux.fr\\
\url{https://www.math.u-bordeaux.fr/~alasjaun/}
\\
\end{tabular}


\newpage
\begin{thebibliography}{99}

\bibitem{Euler} Leonhard Euler, \emph{De Fractionibus Continuis Dissertatio (1737)}
\newline Translated by M. F. Wyman and B. F. Wyman, \emph{An Essay on Continued Fractions}, Mathematical Systems Theory \textbf{%
18} (1985), 295--328.

\bibitem{OP} Oskar Perron, \emph{Die Lehre von der Kettenbr\"{u}chen (1913)}, Chelsea Publishing Company, New York (last edition 1950).

\bibitem{HW}G. H. Hardy and E. M. Wright, \emph{An Introduction to the Theory of Numbers (1938)}. Oxford University Press (last edition, 2008).

\bibitem{H} H. Hasse, \emph{Zahlentheorie (1949)}, Akademie-Verlag, Berlin (last edition 1963).

\bibitem{BS1} L. E. Baum and M. M. Sweet, \emph{Continued fractions of
algebraic power series in characteristic }$2$\emph{,} Ann. of Math. \textbf{%
103} (1976), 593--610.

\bibitem{BS2} L. E. Baum and M. M. Sweet, \emph{Badly approximable power
series in characteristic $2$}, Ann. of Math. \textbf{105} (1977), 573--580.

\bibitem{MR} W. Mills and D. P. Robbins, \emph{Continued fractions for
certain algebraic power series}, J. Number Theory \textbf{23} (1986),
388--404.

\bibitem{S} W. Schmidt, \emph{On continued fractions and Diophantine
approximation in power series fields.} Acta Arith. \textbf{95} (2000),
139--166.

\bibitem{L} A. Lasjaunias, \emph{A survey of Diophantine approximation in
fields of power series.} Monatsh. Math. \textbf{130} (2000), 211--229.
   
\bibitem{T} D. Thakur, \emph{Function Field Arithmetic.} World Scientific
(2004).

\bibitem{LY} A. Lasjaunias and J.-Y. Yao, \emph{Hyperquadratic continued
fractions in odd characteristic with partial quotients of degree one}, J.
Number Theory \textbf{149} (2015), 259--284.


\bibitem{L1} A. Lasjaunias, \emph{ A note on hyperquadratic continued fractions in characteristic 2 with partial quotients of degree 1}, Acta  Arith. \textbf{178} (2017), 249--256. 

\bibitem{M} K. Mahler, \emph{On a theorem of Liouville in fields of positive characteristic.} Canadian J. Math. \textbf{1} (1949),
388--404.


\bibitem{R} K. F. Roth, \emph{Rational approximation to algebraic numbers.} Mathematika \textbf{2} (1955), 1--20.

\bibitem{U} S. Uchiyama, \emph{Rational approximation to algebraic functions.} Proc. Japan Academy \textbf{36} (1960), 1--2.

\bibitem{V} J.-F. Voloch, \emph{Diophantine approximation in positive
characteristic.} Period. Math. Hungar. \textbf{19} (1988), 217--225.

\bibitem{dM} B. de Mathan,\emph{ Approximation exponents for algebraic functions}, Acta Arith. {\bf 60} (1992), 359--370.

\bibitem{BL} A. Bluher and A. Lasjaunias,\emph{ Hyperquadratic power
    series of degree four}, Acta Arith. {\bf 124} (2006), 257--268.

\bibitem{LdM} A. Lasjaunias and B. de Mathan, \emph{Thue's Theorem in positive Characteristic.} Journal fur die reine und angewandte Mathematik \textbf{473} (1996), 195--206.


\bibitem{BR} M. Buck and D. Robbins, \emph{The continued fraction expansion
of an algebraic power series satisfying a quartic equation}, J. Number
Theory \textbf{50} (1995), 335--344.

\bibitem{L2} A. Lasjaunias, \emph{Diophantine approximation and continued
fraction expansions of algebraic power series in positive characteristic},
J. Number Theory \textbf{65} (1997), 206--225.


\bibitem{L3} A. Lasjaunias, \emph{Continued fractions for hyperquadratic
power series over a finite field}, Finite Fields Appl. \textbf{14} (2008),
329--350.

\bibitem{L4} A. Lasjaunias, \emph{ On Robbins' example of a continued fraction for a quartic power series over $\F_{13}$}, J. Number Theory {\bf 128} (2008), 1109--1115.

\bibitem{L5} A. Lasjaunias, \emph{ Algebraic continued fractions in $\F_q((T^{-1}))$ and recurrent sequences in $\F_q$ }, Acta  Arith. {\bf 133} (2008), 251--265. 
 
\bibitem{AL} Kh. Ayadi and A. Lasjaunias, \emph{ On a quartic equation and two families of hyperquadratic continued fractions in power series fields}, Moscow Journal of Combinatorics and Number Theory {\bf 6} (2016), 132--155. 
 
\end{thebibliography}
\end{document}